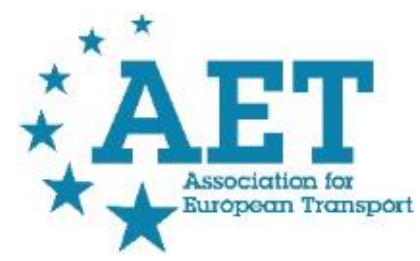

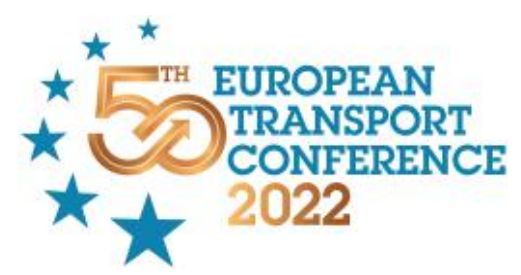

# A METHODOLOGICAL APPROACH FOR DESIGN OF ENERGY EFFICIENT URBAN LOGISTICS SYSTEMS

Emilienne Lardy, Mariam Lafkihi, Eric Ballot
Mines Paris, PSL University, Centre for management science (CGS),
i3 UMR CNRS, 75006 Paris, France

## 1. INTRODUCTION

One of the main challenges for the future in logistics lies in energy consumption. France - as many other countries - has the ambition to reach net zero carbon emissions by 2050. To this aim, the French Electricity Transportation Network has published the most ambitious prospective analysis to date *Energy Pathways 2050,* assessing multiple scenarios involving various energy mixes (wind, solar, nuclear etc.) (Réseau de Transport d'Électricité, 2021). It concludes that the net carbon neutrality in 2050 can only be made possible through energy sobriety, despite mobility switch to electricity (i.e. at least a 40% decrease final energy consumption, plummeting from 1600 TWh to 930 TWh in 2050).

In France, transportation is the second most energy-guzzling field, between the residential and the industrial demands (Ministère de la Transition Ecologique, 2021). While regional and national freight transport bet on the densification and massification of flows and on high capacity vehicles with a view towards energy-efficiency, an opposite trend has been highlighted in urban areas. The city logistics' objective of reducing both public disturbance and atmospheric pollutants (including $CO_2$) entails two simultaneous tendencies. Firstly, the shift in fuels, swapping gas for biofuels, electricity, and even human propulsion. Secondly, the reduction of vehicle capacity, increasing the use of light commercial vehicles for the last mile trips (i.e. Gross Vehicle Weight <3.5 tonnes) (Oliveira et al., 2017). The paradox lies in the abovementioned fact that the shift in fuel mix may not be able to cope with increased overall energy demand resulting from smaller vehicles, more trips and more distance covered.

In view of the reasons discussed above, this work does not seek to quantify energy consumption in city logistics as many other works, but asks how look upon urban freight distribution if energy sobriety is bumped as its top priority, under the specific constraints of dense areas. Instead of confining to the optimization of current practices as a mean of improvement, this work seeks to open widely the range of design parameters explored as part of the urban distribution scheme. This can involve taking advantage of new technologies and innovative processes as well as putting back on the agenda some discontinued practices.
While early stage design offers the biggest leeway to influence impacts (Midler, 1991), the most established impact evaluation methodologies require precise knowledge of the system as well as its components and materials (e.g. Life-Cycle Analysis, Energy Process Analysis). This dilemma has been recognized in many fields, including civil engineering, motivating macro-component based and simulation based approaches (Gervásio et al., 2014), as well as in many analytic approaches (Carlos and Nepomuceno, 2012). We side with the latter and we provide a new methodology devised to comprehensively assess diverse scenarios for freight transportation through the prism of energy sobriety and transportation its major energy component.
In this work we show that taking a step back to fundamental energy theorems yields an analytical formulation of mechanical energy losses from freight distribution. It includes transportation modes characteristics (e.g. aerodynamics, rolling resistance etc.), and kinetics repercussions of delivery management, and infrastructure accessibility (e.g. driving cycles, congestion etc.).
These features make our methodology apt to consider and to contrast diverse worthwhile urban freight distribution schemes. Namely, we chose to examine the exploitation of a global dedicated scheduled system, as an alternative to scattered road haulage subject to traffic and congestion. To this end, we investigate synergies in freight consolidation and intermodal transport either aboveground or underground. Indeed, an intermodal system cannot be limited to a mere transportation mean because, due its interfaces between transportation modes that can incorporate several loading and unloading stations, it can also be considered as an integrated sorting system. Thus it can lend itself well to the regrouping of freight and to the interconnection of supply chains in the urban area. It has been proven in maritime transport industry that standardized containers can enable seamless transshipment in intermodal transport. The same could apply to urban logistics by means of urban containers.

Our contribution is twofold. Firstly, we propose an analytical approach to characterize fundamental energy losses from urban freight distribution suited for early stage design. Secondly we put in practice the approach to compare current practices with an alternative scheme consisting in a global intermodal distribution system, and to draw conclusions regarding best practices by means of an impact assessment of the main design parameters, highlighting the potentials.

## 2. LITTERATURE REVIEW

As soon as the 70's, the cost of energy urged researchers to dive into comparative studies to quantify energetic ramifications of logistics choices. Fels (1975) was among the first who compared the energy consumed per passenger-mile and per seat-mile for several urban passenger transportation modes, including operational consumption as well as embodied energy. He found that rapid rails and city buses are the least energy-guzzling on both criteria, apart from bicycles and walk. So far, this branch of works has been updated to account for the changing circumstances. Thus, McKinnon (2011) listed both the observed progress and the remaining perspectives for improvement on vehicle utilization and energy efficiency of transport. Meanwhile, Figueroa et al. (2014) broke down the energy consumption among transportation modes and notes that the said progress tends to be offset by the switch to more energy intensive transportation modes (e.g. rail and ships to trucks and planes). For these works the typical way of deriving the total energy consists in multiplying a relevant variable (e.g. miles, tonnes, number of vehicles) by a unit consumption coefficient. It is either expressed directly (e.g. kJ, MWh) or through the corresponding amount of fuel (e.g. L of petrol, kg of LNG). The said coefficient is an averaged value that hides disparate instances (e.g. speed, payload, slope…), that can be put into light by detailed energy models based on fundamental mechanics.
Such models have been developed and used for targeted parameters' measurements and validations, as well as to assess the performances of specific technologies applied to transportation. In this way, Hunt et al. (2011) proposes and validates a model for fuel consumption of heavy vehicles, which is built upon by Odhams et al. (2010) by means of an impact assessment of the parameters. Andriaminahy et al. (2019) compares three estimation methods for aerodynamic drag and rolling resistance of light vehicles. These authors prove their results to be accurate with a 1% error margin for laden vehicles (7.7% for unladen vehicles in Hunt et al.). Meanwhile, Van Sterkenburg et al. (2011) examines regenerative braking efficiency of an electric public transport vehicle and an electric garbage truck. Guandalini and Campanari (2018) compares battery and hydrogen fuel cell for light and heavy electric freight vehicles. Stolaroff et al. (2018) audits the use of drones for package delivery. In short, fundamental mechanics have mainly been a tool to refine model for specific technologies. But the recent work from Ankur et al. (2022), on the contrary, leverages these models across various transportation modes (e.g. passenger electric buses, trams, trains, coaches and freight electric semi-trucks, trains, flights) as a mean to generate a versatile model for fuel consumption (i.e. applicable to all transportation modes) in the same vein as our endeavour.
Concurrently, works have been set up to assess the performances of courses of action. In this way, driving cycles have been utilized to describe driver's behaviours and their impact on fuel consumption. In Europe, André has collected 70000 kilometres of real world data (André, 1991) and has drawn out representative driving cycles (André, 2004). These driving cycles - coined 'Artemis' - and their counterparts from all around the world and from various types of vehicles have been catalogued by Barlow et al. (2009). In an effort to make sense of the plethora of driving cycles - including official homologation cycles and non-official real world ones - Zaccardi and Le Berr (2013) use data analysis techniques to identify the few most characteristic ones. This branch of research is still going strong in hopes of refining the quality of representativeness of driving cycles; Giraldo et al. (2021) puts forwards indicators of said representativeness (e.g. Relative Differences, Mean Relative Absolute Difference) and concludes that a minimum driving cycle duration of 25 minutes is essential. On analogous grounds, Boggio-Marzet et al. (2021) examine the delivery route typology's repercussions on energy efficiency, by means of data collection in Madrid coupled with an analytic model. They find that the city centre routes - in comparison with peri-urban areas - are plentiful with traffic lights, pedestrian crossings and congestion, which fallout include low speeds and low energy efficiency.
Our work subsumes these research branches into a comprehensive versatile model for urban freight distribution. We focus on the mechanical energy losses from transportation modes suited for urban areas, adding the option of going through a tunnel - as can be necessary in high land-pressure areas, and we layer the resulting formulation through all the distribution steps from outskirts to delivery.

## 3. METHODOLOGY

The methodology used in this work is devised to assess the engine energy supply required for the core function of moving a freight unit load from the outskirts of the city to its destination. All the available implementation choices - both in technology and in layout - have to be translatable in our methodology. The method consists in taking a step back to fundamental energy theorems, which parameters can be set to stand for any scenario of interest. The methodology (Figure 2) includes transportation modes characteristics (e.g. empty and loaded weights, aerodynamics etc.), geographical features (e.g. distances, slop etc.) as well as kinetics repercussions of delivery management, and infrastructure accessibility (e.g. driving cycles, congestion, load rate etc.). These specifications yield, by means of the work-energy theorem, the theoretical energy to be provided to compensate energy losses.

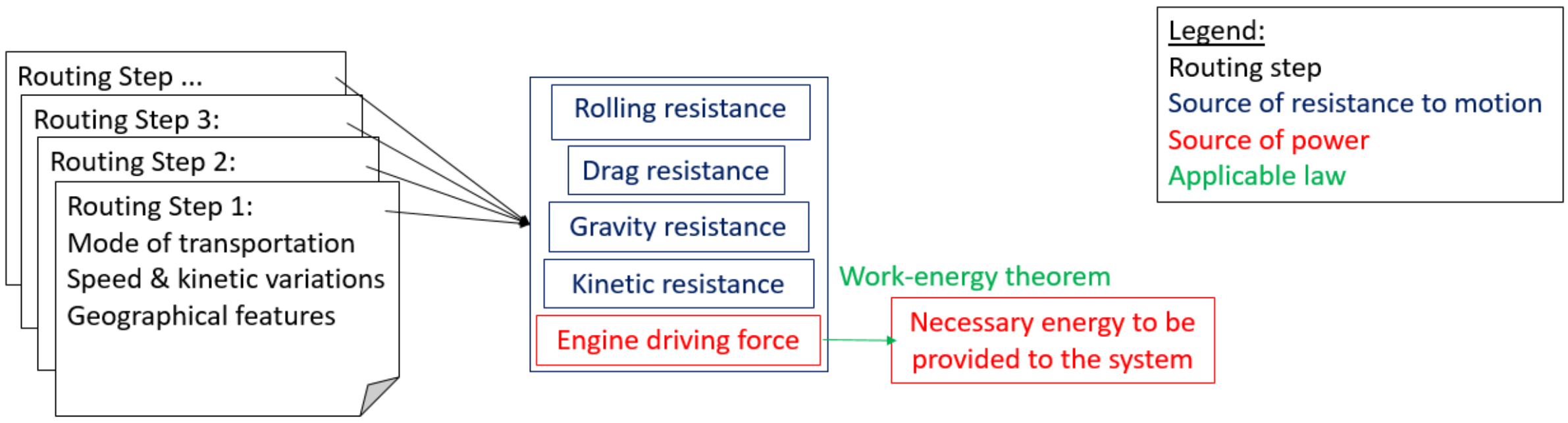

*Figure 1: Energy loss modelling from fundamental energy theorems*

### 3.1 Energy loss modelling

The sources of resistance to motion to be taken into account as part of mechanical energy losses are as follows:

| Type | Formulation | Notation |
|---|---|---|
| Rolling resistance | $F_r = C_r * m * g$ | $C_r$ *rolling resistance coefficient of the vehic*<br>$m$ *total mass of the vehicle*<br>$g$ *acceleration of gravity* |
| Aerodynamic drag *in open air* | $F_a = \frac{1}{2} * \rho * v^2 * S * C_a$ | $S$ *frontal surface of the vehicle*<br>$C_a$ *aerodynamic coefficient of the vehicle*<br>$\rho$ *volume weight of air*<br>$v$ *speed* |
| Aerodynamic drag *in a tunnel* (Vardy, 1996) | $F_a = \frac{1}{2} * \rho * v^2 * \left[A_t * (k_N + k_T) + f_z * l_z * L_z * \left(1 + \frac{At}{Ann}\right)\right]$ | $k_N$; $k_T$ *nose and tail loss coefficient*<br>$f_z$ *skin friction*<br>$l_z$ *cross sectional perimeter of trolley*<br>$L_z$ *lenth of trolley*<br>$A_t$ *cross sectional area of tunnel*<br>$A_{nn}$ *cross sectional area of annulus*<br>$\rho$ *volume weight of air*<br>$v$ *speed* |
| Gravity resistance | $F_g = m * g * \sin(\alpha)$ | $m$ *total mass of the vehicle*<br>$g$ *acceleration of gravity*<br>$\alpha$ *angle between the road and the horizontal* |
| Kinetic resistance | $F_k = m * a$ | $m$ *total mass of the vehicle*<br>$a$ *accleration of the vehicle* |
| | It is approached by driving cycles' Art.Kinema parameters theorized by (Haan and Keller, 2004) and catalogued by (Barlow et al., 2009). The total positive kinetic energy is: | |
| | $PKE = \frac{1}{dist} * \sum i \begin{cases} v_i^2 - v_{i-1}^2 \ if\ (v_i > v_{i-1}) \\ 0\ (else) \end{cases}$ | *dist total distance*<br>$(v_i)$ *itemized speeds* |

### 3.2 Delivery route modelling

Our work applies the previously described methodology to compare current practices with an alternative energy aware scheme, and to pinpoint the differences in their outcomes in terms of energy loss.

Essentially, the alternative scheme we examine consists in a two-stage network ( Figure 2 ), involving consolidation centres on the outskirts of the city, and proximity transshipment points inside the city. The former and the latter are connected by a railway infrastructure that can either be aboveground or underground when the land pressure is high. The urban area is paved into catchment areas, which centroid is the ideal location for the proximity transshipment -and contingently resurfacing- point. The distribution system is rounded off by short delivery trips with a 'soft' mode of transport, from the nearest proximity transshipment point to the destination.

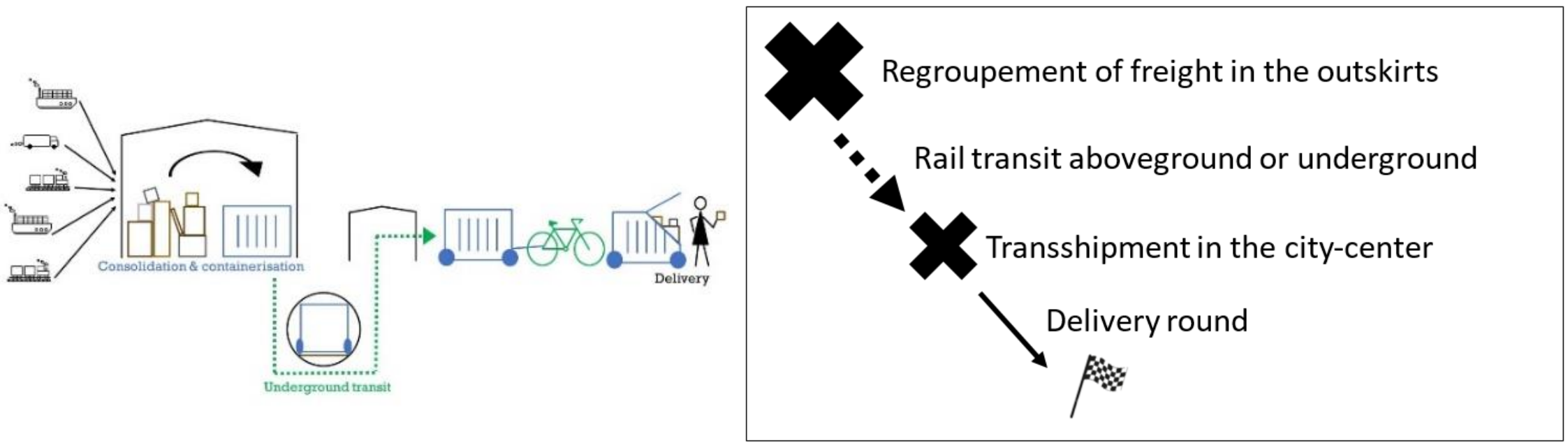

*Figure 2: Illustration of a dedicated urban freight supply system involving consolidation on containerization*

For purposes of generality, the urban layout is simplified as a general dense area accessed via a ring road. This allows for two distinct set of routing steps corresponding to the current practices versus the alternative scheme, illustrated on Figure 3. It is applied to the case of Paris for quantifying all relevant figures of freight supply, benefiting from previous studies and surveys, as well as government approved open data.

In the "current practices" scheme, freight is shipped from the departure point (i.e. distribution centre), to its destination via the ring road, entering and exiting the ring road through the closest ring road doors. The vehicle mix is, as found by the survey 'Marchandises en Ville' by the LAET (Serouge et al., 2018), approximately 30% van, 30% light commercial vehicles, 30% rigid truck, 10% articulated vehicles. The standard delivery rounds include 35 stops, scattered on an average area of two Paris districts (i.e. out of a total of 20 districts) - these figures were procured by logistics providers.

In the "alternative scheme", from the same departure point, a detour is made to reach the closest entry point of the dedicated railway, exploiting essentially heavy commercial vehicles as they will not enter the urban area. On site, the freight is regrouped, scheduled, and transshipped to a rail trolley, routed through the railway network to the nearest exit to its destination point. This step includes not only transportation but also the lift down to and up from the underground transport level - considering that at least a section is indeed underground. The last trip is carried out by light commercial vehicles or cargo bikes.

The distances are based on geometric calculation for steps 1, 2 and 3. For the delivery round -step 4-, we rely on the continuous approximation model provided by Daganzo (Daganzo, 1984) to estimate the total delivery distance:

$$D = 2 * e * \frac{N}{C} + k * \sqrt{N * S} \text{ with } \begin{cases} N \text{ number of destinations in the catchment area} \\ C \text{ number of stops per vehicle} \\ e \text{ expected distance to reach the destinations} \\ S \text{ catchment area} \end{cases}$$

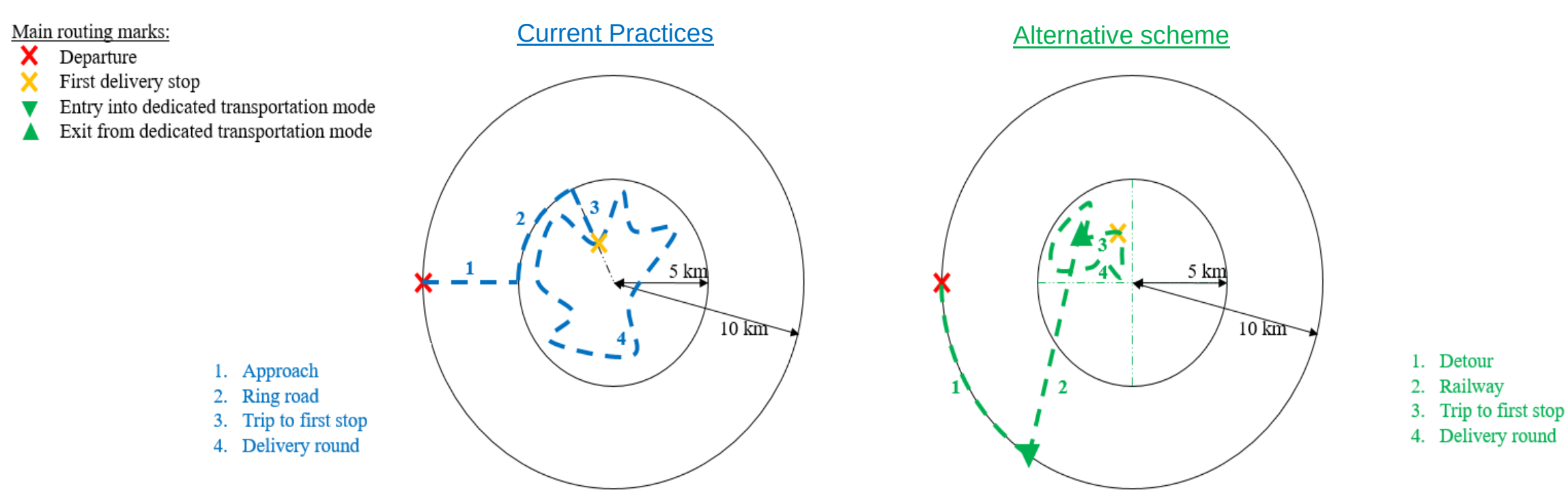

*Figure 3: Main routing steps for freight supply into Paris*

## 4. ANALYSIS AND RESULTS

The methodology presented beforehand establishes an analytic model of energy losses due to transit and delivery of freight. It is applied to two contrasting schemes: the current practices for freight supply in Paris, and a global intermodal system. The methodology's subjacent layers bring into light both the origin and the nature of energy losses. The contrast between both schemes on these matters will be discussed here and put into perspective as part of current city logistics' developments.

### 4.1 Technology related learnings

In terms of energy sobriety among freight land transportation modes for long distances, trains win hands down (substantiated recently by Ankur et al. (2022)). The reason is that aerodynamic and rolling resistances are minimal next to the inertia of the motion of such a large mass. This allows to approximate the energy consumption as a function of the speed and the number of stops - the distance covered disappears. European Cooperation in the Field of Scientific and Technical Research (Organization) and European Commission, 1999, p. 87 gives

$$E = k_1 * \frac{V^2}{\ln(x)} + k_2 \text{ where } \begin{cases} E \text{ is the energy consumption in } kJ \\ v \text{ is the average speed of the train} \\ x \text{ is the distance between stops in } km \\ k_1 \text{ and } k_2 \text{ are constants} \end{cases}.$$

But large trains are not suitable for urban freight distribution, they would have to be substituted - for flexibility purposes - by small capacity rail trolleys. The hegemony of railways for moving large masses over long distances thus disappears, and it is yet to be proven if railways hold their appeal, energy-wise, in the context of small rail trolleys over short urban distances.
To that aim, we compare energy losses of currently used heavy and light freight trucks with rail trolleys of equivalent payload. To provide a comparison baseline, we will consider the translational motion of 1 tonne of freight over 1 kilometre, including the start-up and stop motions.

In turns out rail trolleys can curb the explosive increase in energy demand from the flexibility and the seamlessness provided by smaller vehicles. Indeed, while high capacity rail and road vehicles have similar energy outputs, the rail trolleys' performances (solid lines in Figure 4) seem to deteriorate very slowly with capacity reduction compared to road vehicles' performances (dotted lines).

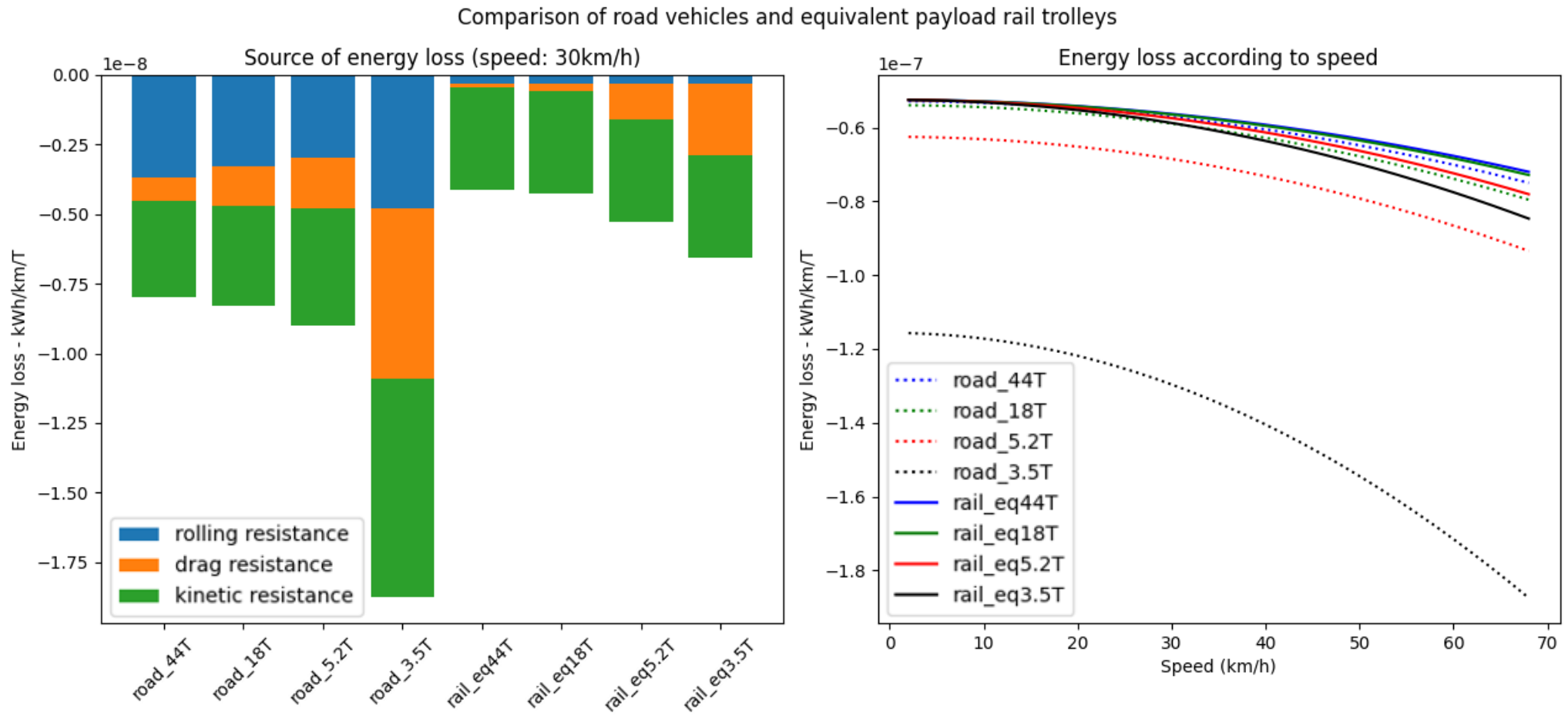

*Figure 4: Comparison of energy losses from road vehicles and equivalent payload rail trolleys*

The breakdown of energy losses shows that rail trolleys have very low rolling resistance, consequently the aerodynamic drag is prevalent. In considering the exploitation of small rail trolleys, the urban context where speeds are low -in Paris they currently range between 15km/h and 30km/h- and the aerodynamic drag is weak, is particularly worthwhile.

## 4.2 Traffic flow related learnings

The decomposition of energy losses both between routing steps, and between sources of resistance is presented in Figure 5, for worst case traffic and best case traffic.
It turns out that kinetic energy losses are strikingly predominant, based on initial assumptions they account for 98% of total direct energy losses in current practices and for 97% in the alternative scheme. The alternative scheme profits from a dedicated transportation mode, which bypasses congestion in the ringroad and the entry into the city, and allows a constant speed throughout. It induces between a 55% cut and a 67% cut in energy losses. The alternative is also more resilient to bad traffic conditions as the total energy losses are increased by 32% between best and worst case traffic, against 50% for current practices. These stark figures are robust to small inexactness. Still, we are working on replacing the use of driving cycles with real data on speed and acceleration in Paris to get more precise figures.

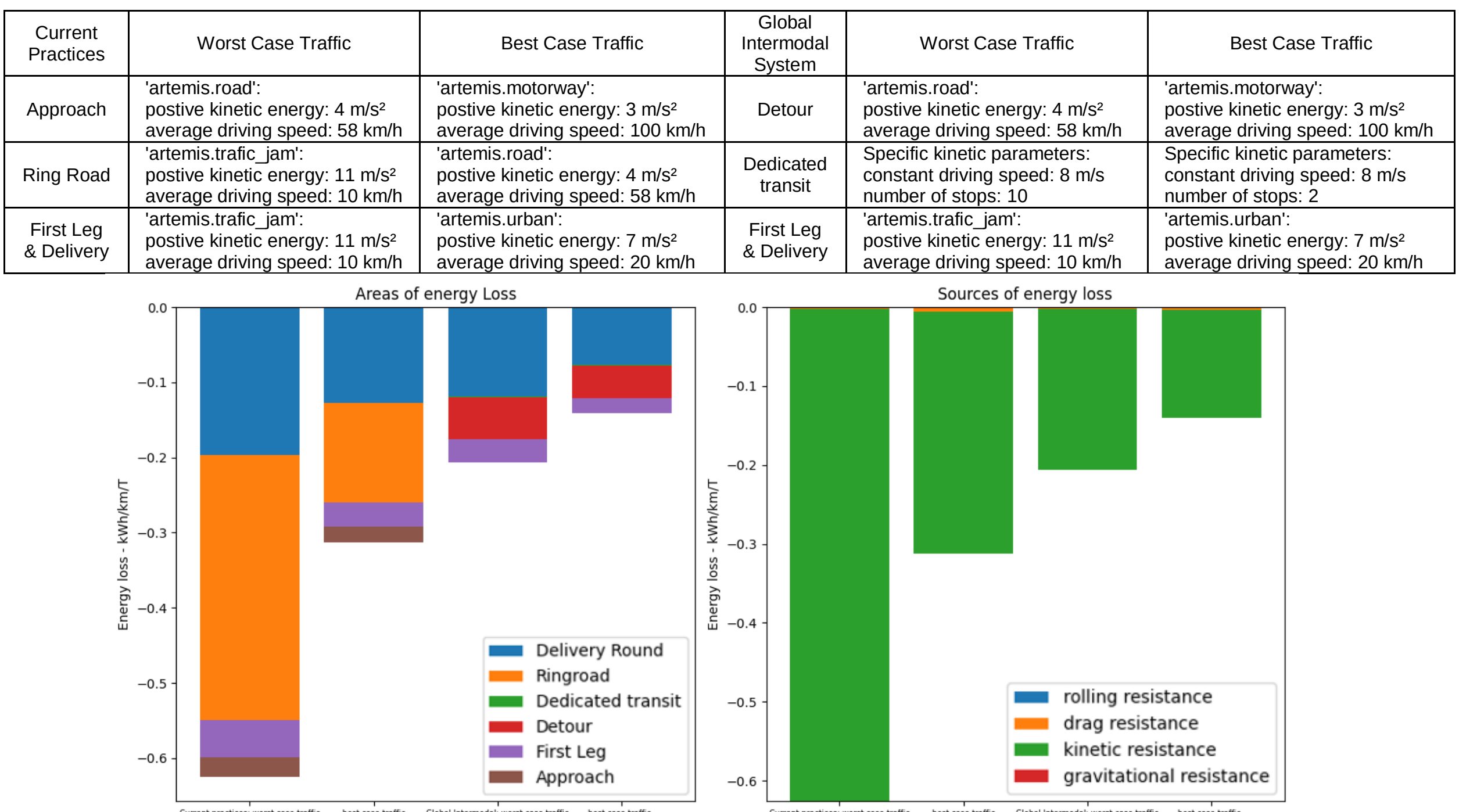

| Current Practices | Worst Case Traffic | Best Case Traffic | Global Intermodal System | Worst Case Traffic | Best Case Traffic |
|---|---|---|---|---|---|
| Approach | 'artemis.road':<br>postive kinetic energy: 4 m/s²<br>average driving speed: 58 km/h | 'artemis.motorway':<br>postive kinetic energy: 3 m/s²<br>average driving speed: 100 km/h | Detour | 'artemis.road':<br>postive kinetic energy: 4 m/s²<br>average driving speed: 58 km/h | 'artemis.motorway':<br>postive kinetic energy: 3 m/s²<br>average driving speed: 100 km/h |
| Ring Road | 'artemis.trafic_jam':<br>postive kinetic energy: 11 m/s²<br>average driving speed: 10 km/h | 'artemis.road':<br>postive kinetic energy: 4 m/s²<br>average driving speed: 58 km/h | Dedicated transit | Specific kinetic parameters:<br>constant driving speed: 8 m/s<br>number of stops: 10 | Specific kinetic parameters:<br>constant driving speed: 8 m/s<br>number of stops: 2 |
| First Leg & Delivery | 'artemis.trafic_jam':<br>postive kinetic energy: 11 m/s²<br>average driving speed: 10 km/h | 'artemis.urban':<br>postive kinetic energy: 7 m/s²<br>average driving speed: 20 km/h | First Leg & Delivery | 'artemis.trafic_jam':<br>postive kinetic energy: 11 m/s²<br>average driving speed: 10 km/h | 'artemis.urban':<br>postive kinetic energy: 7 m/s²<br>average driving speed: 20 km/h |

*Figure 5: Breakdown of energy losses for current practices and the examined alternative scheme*

From the predominance of kinetic losses, we can deduct that there is no better line for action, among initiatives with a plausible short term and large scale implementation, then avoiding traffic congestion. This immediately echoes the strategy of off-peak hours' deliveries (OPHD), in respect of which Sánchez-Díaz et al. (2017) provide a review of theory and practice. The concept is to shift the delivery trips during the night, in order to avoid peak hours' traffic congestion, and to save both travel time and fuel, as well as diminishing public disturbance. Without matching the performances of a dedicated scheduled transportation mode (scattered vehicles cannot be cadenced, plus crossroads and traffic lights remain and generate speed variations), it seems that OPHD has high potential in terms of energy sobriety. The above-mentioned review reports that the travel time saved can exceed 50% according to pilots in London, Paris and Stockholm, exhibiting the aptitude to uphold speed and to cut kinetic losses.

## 4.3 Spatial design related learnings

From the precedent section we can predict that the distance efficiency of the railway network has very little impact on energy losses. By distance efficiency, we refer to the ratio between the real distance allowed by the network between the entry and the exit points, and their distance as the crows fly (corresponding to a web-like network in which there is a direct and straight connection between each entry point and each exit point). Because of the negligible kinetic energy losses of a dedicated and cadenced railway, any distance surplus inside the railway network has a marginal impact on total energy losses. As an illustration, Figure 6 compares the total effect indices of variance-based sensitivity analysis of the distance efficiency (D), the number of access points (A) and the number of exit points (E), over both the distance covered and the energy losses. The graphs are displayed using the SALib library (Herman and Usher, 2017).

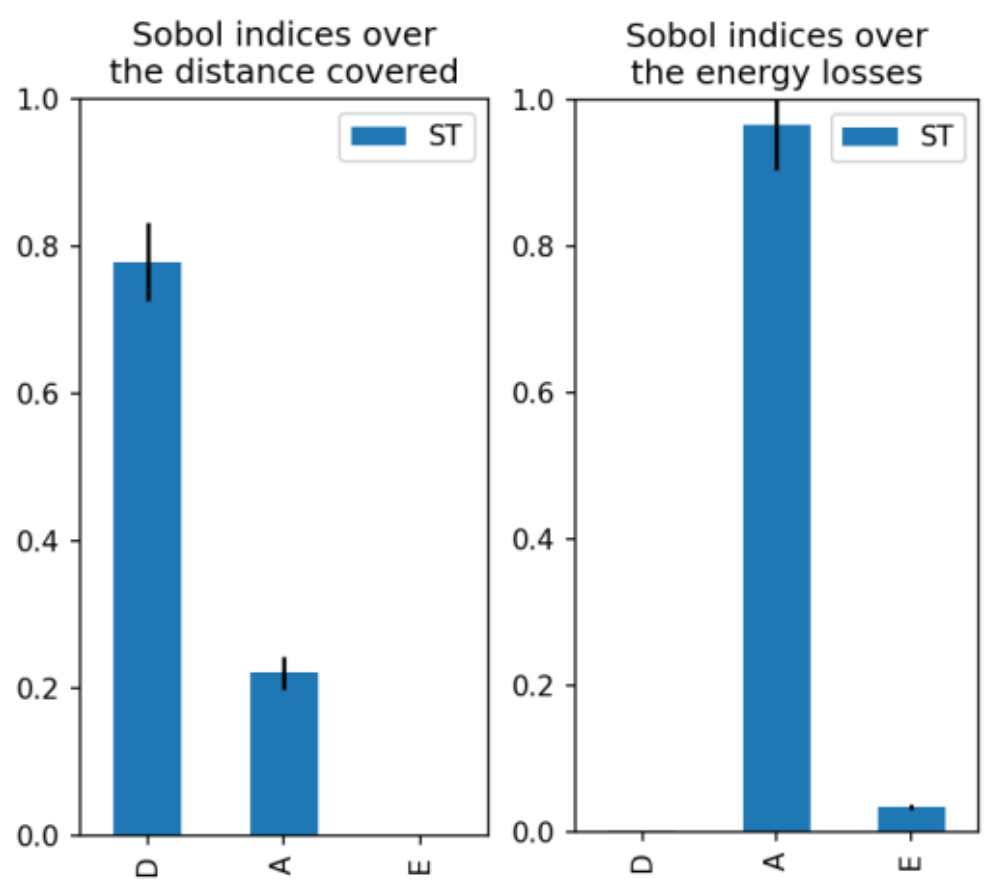

*Figure 6: Sobol Indices of the distance efficiency (D), the number of access (A) and exit (E) points*

This is of crucial importance as it means that minimizing the number of connections between the nodes of the network - provided that its throughput can absorb the flow - will yield low embodied energy without degrading much the direct energy losses.

### 4.4 Management practices related learnings

The perspectives from collaboration can be delved into through two meaningful parameters which are expected to benefit from pooled deliveries, namely the density of destinations and the average load rate. The default values describing current practices are, as mentioned in 3.2, a 3.3 density (i.e. 35 delivery destinations per two Paris districts), and a 40% average load rate (i.e. 70% on headhaul and 10% on backhaul).
Figure 7 shows the repercussions on both schemes of these two parameters. It confirms that the global intermodal scheme outperforms the current practices regardless of the parameters' values; this is natural as the alternative scheme cuts the distance covered by high intensity vehicles. The intermodal scheme seems more sensitive, as a 10% improvement on each parameter yields a 19% contraction in energy losses, against only 15% for current practices. This gap is thin, but while parameters' improvement is within reach in both case - if anything through routing optimization - the range of improvement is far greater in a global system, as expected from works on horizontal collaboration in logistics (Shenle Pan et al., 2019). For a global urban distribution system, no feedback is available to give tried and validated figures for the parameters, as collaboration in urban distribution has never been achieved on a wide scale. Upper bounds could be set as high as a 53% load rate (i.e. 90% on headhaul and 15% on backhaul) and a 55#/km² density. This density is based on the 93725 sales point identified by the Paris open data web portal (accessed in April 2021) out of a 104.5km² area, and considering 1h delivery windows throughout a 16h time range. As at happens, while consolidation of city supply flows knows stark barriers to implementation (Björklund and Johansson, 2018), it has been proven that associating a dedicated transit makes it worthwhile for a wide range of good (Boerkamps et al., 2000), making the global intermodal scheme coherent and plausible.

As for higher average load rates, they can only steam from large scale valuation of backhauls, which is the case in a circular economy but which goes beyond the scope of logisticians' own management practices.

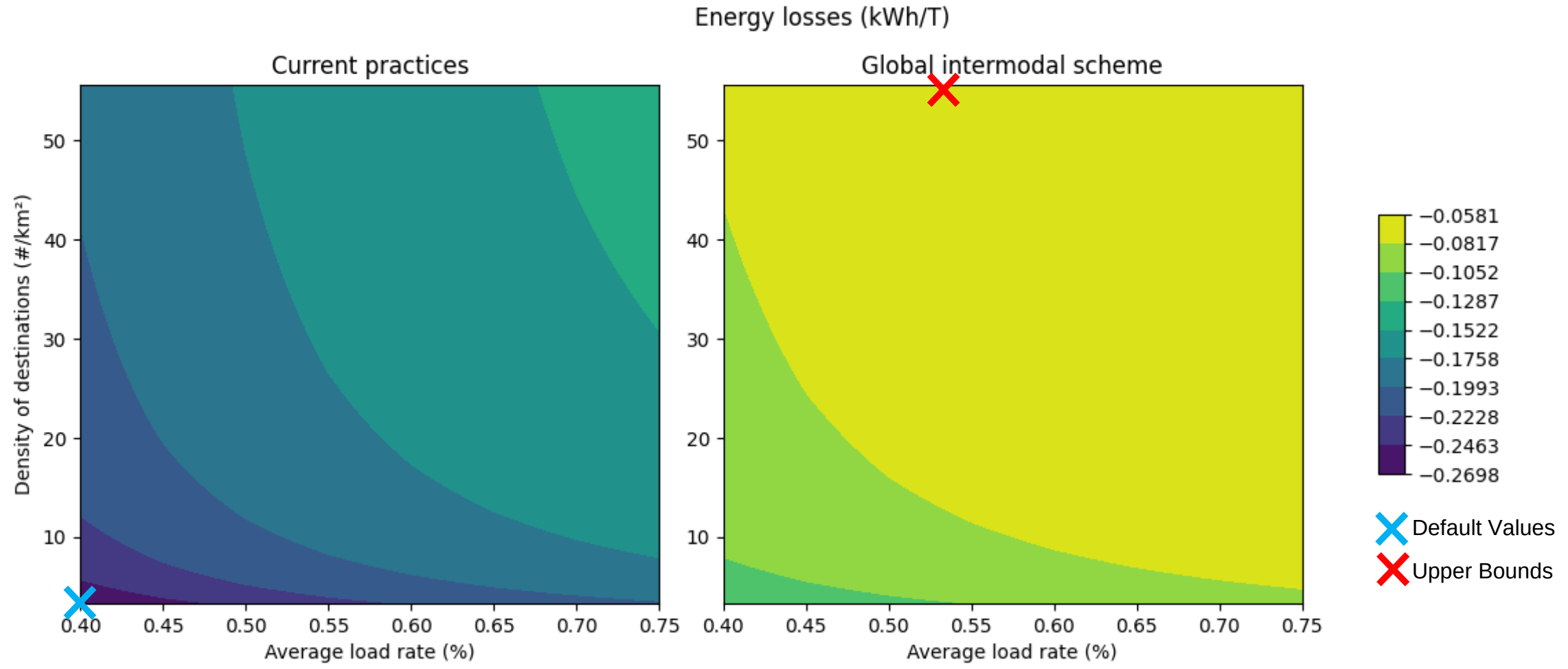

*Figure 7: Repercussions of the average load rate and the density of destinations*

## 5. DISCUSSION

The first point to address is the distinction that we draw between mechanical energy losses and energy consumption, the former leading to the latter by multiplication with the drivetrain efficiencies and addition with the auxiliary energy and the HVAC energy (Ankur et al., 2022). Our works sticks to mechanical energy losses for two reasons. Firstly, raw mechanical energy loss is eloquent to depict the ins and outs of energy sobriety. Secondly it avoids speculating on technical progress' transferability, as the current improvements in drivetrain efficiencies hold no more to transportation modes themselves than to their markets turnover and subsequent investments in research and development.

The second point to address is the blindness of our methodology on embodied energy, which constitutes a considerable limit. Mainly, the energy invested in the construction of the equipment and infrastructure is considerable in the overall energy cost of a transportation system (Levinson et al., 1984). To stay coherent with our goal of energy sobriety, we consider it a priority to revalue currently underused urban facilities (unused railways, warehouses, underground parking in Paris or other). The necessary infrastructure to the examined logistics system should be embedded in the existing urban surroundings, and entail as little new constructions and extra artificialization of soils as possible.

This applies, for example, to the global intermodal system for urban distribution that we have been examining. As it happens, an open distribution system as the one examined in this work frees the logistics providers of the last mile dispatching and scheduling by providing a trustworthy alternative. It therefore curbs their need for distribution centres in the outskirts of the city, resulting in unused warehouse area - including storing and sortation systems. It can even be expected that the shift to a global system leads to less total required resources for the separation and transshipment process, which is a known benefit of collaboration (Ballot et al., 2014). This is why the necessary amount of civil engineering cannot be assimilated to -and is nowhere near- the bottom up construction of system's infrastructure. Instead, a dedicated study should discriminate the existing infrastructure that can be redirected towards making up the global supply system, from the missing infrastructure that actually has to be built.
In order to include embodied energy in our methodology, this thought process has to be converted into a general and replicable methodology. This work is out of the scope of the research project presented here.

## 6. CONCLUSION

In the objective of providing a straightforward overview of the origin and the nature of energy losses for broad range of early-stage urban freight distribution schemes, we have proposed a methodology. It subsumes the transportation modes characteristics, the geographical features and the kinetics repercussions of distribution management into an analytic formulation by means of the work-energy theorem. This makes for an easily accessible tool befitting impact assessments and comparative studies. The methodology was put into practice via application to two contrasting schemes: current practices for freight distribution and a global intermodal system. This application was ground for exploring the range of inquiries and learnings within reach due to the subjacent layers of the methodology.

## 7. CONCLUSION

The authors thank GEODIS for funding and collaborating with the Physical Internet Chair, as well as and VINCI Construction, VINCI Concession and Île de France for their support during the AMI program "Fret et Logistique".

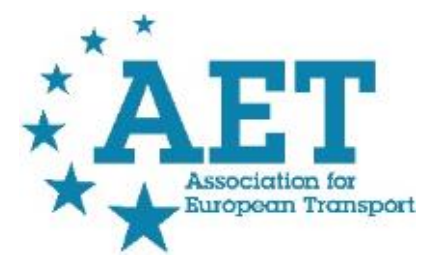

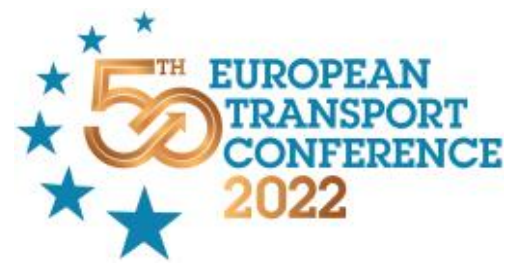

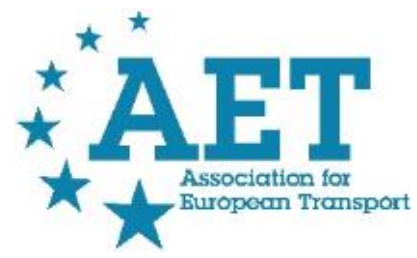

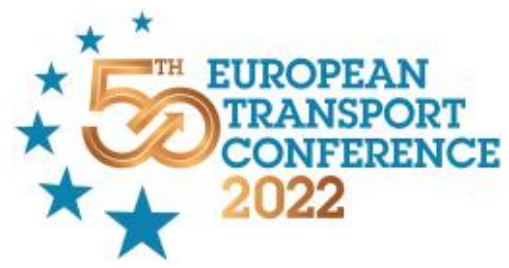